%% file: Matula2.tex
\documentclass[english]{amsart}
\usepackage{babel}
\usepackage[latin1]{inputenc}
\usepackage[pdftex]{color,graphicx}
\usepackage{amssymb}
\usepackage{exscale}
\usepackage{latexsym}
\usepackage{color}
\usepackage[all]{xy}
\usepackage{fancyhdr}
\usepackage[mathcal]{eucal}
\usepackage{amsmath,amsfonts}
\usepackage{eufrak}
\usepackage{multicol}
\usepackage{epstopdf}
\usepackage[pagewise]{lineno}
\usepackage[headheight=12pt,a4paper,scale=.8,ignorefoot=true,verbose,%
hmarginratio=1:1,%
vmarginratio=1:1]{geometry}
\usepackage{mathrsfs}
\usepackage{enumerate}
\usepackage{longtable}
\usepackage{booktabs}
\usepackage{calc}
\usepackage{pgf, pgfarrows, pgfnodes}
\usepackage{lineno}
\usepackage{fancyhdr}
\newtheorem{prop}{\rm \textbf{Proposition}}[section]

\title{Rooted tree statistics from Matula Numbers}
\author{Emeric Deutsch}
\date{}

\begin{document}
\maketitle
\vspace{-.2in}
\begin{center}Polytechnic Institute of New York University\\emericdeutsch@msn.com\end{center}
\begin{abstract}
There is a one-to-one correspondence between natural numbers and rooted trees; the number is called the Matula number of the rooted tree. We show how a large number of properties of trees can be obtained directly from the corresponding Matula number. 
\end{abstract}

\section{Introduction}

Let \(p_m\) denote the \(m\)-th prime number (\(p_1=2, p_2 = 3, p_3=5, ...\)). We call \(m\) the \emph{order} of \(p_m\). For a positive integer \(n\), we denote by \(\Omega(n)\) the number of prime divisors of n,  counted with multiplicities. 
\\

For a rooted tree \(T\), its \emph{Matula number} \(\mu(T)\) is defined recursively in the following manner \cite{Matula}. For the 1-vertex tree \(T\) we define \(\mu(T)=1\). Otherwise, let \(T_1, T_2, ..., T_d\) be the subtrees of \(T\) rooted at the vertices \(v_1, v_2, ...,v_d\), adjacent to the root of \(T\) (see Fig. 1). Then, we define 
\[
\mu(T)=p_{\mu(T_1)}p_{\mu(T_2)}...p_{\mu(T_d)}.
\]
\begin{figure}
 \begin{center}
  \input{fig1.tex}
 \end{center}
\caption{}
\end{figure}

Conversely, given a natural number \(n\), the above procedure can be done in reverse, leading to a unique rooted tree \(\tau(n)\) with Matula number \(n\). We illustrate this on the arbitrarily chosen \(n=987654321\). We label the root of the desired \(\tau(n)\) by \(987654321=3*3*17*17*379721\). Since we have 5 prime factors, the degree of the root is 5. We draw 5 edges from the root and we label their endpoints, for example from left to right, by the orders of these 5 primes, i.e. by 2, 2, 7, 7, 32277 since \(p_2=3, p_7=17\), and \(p_{32277}=379721\). More conveniently, instead of 32277 we use the label \(32277=3*7*29*53\). Now we repeat this procedure for each of the 5 newly labeled vertices. For example, considering one of the vertices labeled 2, we draw one edge whose endpoint will have the label 1 (the order of the prime 2). This is a leaf of the required tree \(\tau(n)\) and the procedure stops along this branch. The reader is asked to follow in detail the way the tree \(\tau(987654321)\) has been obtained (Fig. 2).  
\\

\begin{figure}
 \begin{center}
  \input{fig2.tex}
 \end{center}
\caption{}
\end{figure}

We will use the following slightly modified description of the above defined correspondence. We give the mapping \(\tau\) from natural numbers to rooted trees, i.e. the inverse of the \(\mu\) function. We define \(\tau(1)\) to be the 1-vertex tree. For a prime number \(p_t\) (the \(t\)-th prime), we define \(\tau(p_t)\) recursively as shown in Fig. 3a. If \(n\) is composite, \(n=rs\) (\(r,s \geq 2)\), then \(\tau(rs)\) is defined recursively to be the tree shown in Fig. 3b, i.e. the trees \(\tau(r)\) and \(\tau(s)\) joined at their roots. Clearly, the obtained rooted tree does not depend on the used factorization \(rs\) of \(n\). 
\\

{\bf Remark.} The prime/composite dichotomy has led us to the same recursive construction of the rooted trees that is mentioned by Czabarka, Sz\'ekely, and Wagner \cite{Czabarka}. Although that paper does not consider the Matula bijection, there may be some overlap when the same statistic is investigated. 
\\

As pointed out by Ivan Gutman and Yeong-Nan Yeh \cite{GutmanYeong}, it is of interest to find statistics on rooted trees \(T\) directly from their Matula numbers \(\mu(T)\). This has been done for several statistics in \cite{ElkGutman}, \cite{GutmanMath}, \cite{GutmanIvic}. 
\\

In this paper we will do the same in a slightly different way and also for several other properties like, for example, various topological indices (Wiener, Zagreb1, Zagreb2, Randi\'c, etc.). We will also consider polynomial-valued statistics (i.e. a finite sequence of statistics), namely the partial Wiener polynomial, the Wiener polynomial, the degree sequence polynomial, and the exit-distance polynomial (the last two terms will be defined below).  
\\

\section{Terminology}

A \emph{rooted tree} is a tree having a distinguished vertex, called the \emph{root}. 
In a rooted tree the \emph{level} of a vertex \(v\) is its distance from the root, i.e. the  length of the unique path from the root to \(v\).
The \emph{height} (called also \emph{depth}) of a rooted tree is the length of the longest path from the root. 
In a rooted tree if vertex \(v\) immediately precedes vertex \(w\) on the path from the root to \(w\), then \(v\) is the \emph{parent} of \(w\) and \(w\) is the \emph{child} of \(v\). 
In a rooted tree vertices having the same parent are called \emph{siblings}. 
A \emph{leaf} in a rooted tree is any vertex having no children. However, in the 1-vertex tree, the root is not considered to be a leaf. Given a rooted tree, the \emph{path length} is the sum of the levels of each of the nodes. The \emph{external path length} is the sum of the levels of each of the leaves (see, for example, \cite{SedgewickFlajolet}). The \emph{distance} between vertices \(i\) and \(j\) of a tree \(T\) is the number of edges on the unique path from \(i\) to \(j\); it is denoted by \(d_T(i,j)\) (or \(d(i,j)\) when no ambiguity is possible). The \emph{diameter} of a tree is the greatest distance between pairs of vertices.  
The \emph{degree} of a vertex \(i\) of a tree \(T\) is the number of edges emanating from \(i\); it is denoted by \(deg_T(i)\) (or \(deg(i)\) when no ambiguity is possible).  
A vertex of degree 1 is called a \emph{pendant vertex} while a vertex of degree \(\geq\) 3 is called a \emph{branching vertex}. A \emph{root subtree} of a rooted tree is any subtree that contains the root. 
\\

The \emph{degree sequence} of a graph is the list of vertex degrees, written in nonincreasing order. This is basically the \emph{comparability index} of a tree introduced by Gutman and Randi\'c \cite{GutmanRandic} (see also \cite{Trinajstic}). For example, for the path on 5 vertices, the degree sequence is \((2,2,2,1,1)\). We define the \emph{degree sequence polynomial} of a graph with vertex set \(\{1,2,...,n\}\) as \(\sum_{i=1}^n x^{deg(i)}\) (it is the generating polynomial of the vertices of the graph with respect to vertex degree). For example, for the path on 5 vertices the degree sequence polynomial is \(2x+3x^2\). 

The \emph{visitation length} of a rooted tree \(T\) is defined as the sum of the number of nodes of \(T\) and of its path length \cite{KeijzerFoster}. 
\\

Given a vertex \(v\) in a rooted tree \(T\), we define the \emph{exit distance} of \(v\), denoted \(\lambda(v)\), to be the distance from \(v\) to the nearest leaf of \(T\) that is a descendant of \(v\). We are interested in three statistics on a rooted tree \(T\), connected with this new concept: (i) the sum of the exit distances of all vertices of \(T\), (ii) the largest exit distance, and (iii) the number of vertices for which the largest exit distance is attained. To obtain all three statistics we define the \emph{exit-distance polynomial} of \(T\) by \(\sum x^{\lambda(v)}\), where the summation extends over all vertices \(v\) of \(T\); it is the generating polynomial of the vertices of the tree with respect to exit distance. To label the vertices of a rooted tree by their exit distances, label the leaves with 0, label their parents with 1, label the so far unlabelled parents of the 1's with 2, label the so far unlabelled parents of the 2's with 3, and so on. Note the monotonicity of the coefficients of this polynomial, due to the fact that each vertex with exit distance \(k\) (\(k\geq 1\)) is the parent of some vertex with exit distance \(k-1\). The reader is asked to verify that the exit-distance polynomial of the rooted tree in Fig. 2 is \(15+9x+5x^2\).      
\\

\begin{figure}
 \begin{center}
  \input{fig3.tex}
 \end{center}
\caption{}
\end{figure}

The \emph{Wiener index} of a tree \(T\) is the sum of the distances between all unordered pairs of vertices of T \cite{Wiener}, \cite{DobryninEntringerGutman}.

The \emph{terminal Wiener index} of a tree \(T\) is the sum of the distances between all unordered pairs of pendant vertices of T \cite{GutmanFurtulaPetrovic}, \cite{SzekelyWangWu}.

The \emph{hyper-Wiener index} of a tree \(T\) has been defined by M. Randi\'c \cite{Randic}, \cite{GutmanChemical}. Later, Klein, Lukovits, and Gutman \cite{KleinLukovitsGutman} have derived the formula \(\frac{1}{2}(\sum d(i,j)^2 +\sum d(i,j))\), where the summations extend over all unordered pairs of vertices. This is now accepted as the definition of the hyper-Wiener index of a connected graph.

The \emph{multiplicative Wiener index} of a tree is the product of the distances between all unordered pairs of vertices of T \cite{GutmanLinertLukovitsTomovic}.

The \emph{Wiener polarity index} of a tree \(T\) is defined as the number of unordered pairs of vertices \({i,j}\) of \(T\) such that \(d(i,j)=3\) \cite{Wiener}.  

The \emph{first Zagreb index} of a tree \(T\) is defined as  \(\sum deg(i)^2\), where the 
summation is over all the vertices \(i\) of \(T\) \cite{GutmanTrinajstic}, \cite{GutmanDasKinkar}, \cite{NikolicKovacevicMilicevicTrinajstic}. 

The \emph{second Zagreb index} of a tree \(T\) 
is defined as  \(\sum deg(i)deg(j)\), where the summation is over all the edges \(ij\) of \(T\) 
\cite{GutmanTrinajstic}, \cite{NikolicKovacevicMilicevicTrinajstic}, \cite{DasKinkarGutman}. 

The \emph{Narumi-Katayama index} of a tree \(T\) is defined as \(\prod deg(i)\), where the product is taken over all the vertices \(i\) of \(T\)
\cite{NarumiKatayama}. 

The \emph{first multiplicative Zagreb index} of a tree \(T\) is defined as \(\prod deg(i)^2 \), where the product is taken over all the vertices \(i\) of \(T\) \cite{GutmanMath}. It is the square of the Narumi-Katayama index. 

The \emph{second multiplicative Zagreb index} of a tree \(T\) is defined as \(\prod deg(i)deg(j) \), where the product is taken over all the edges \(ij\) of \(T\) \cite{GutmanMath}. An other equivalent expression is \(\prod deg(i)^{deg(i)}\), where the product is taken over all vertices \(i\) of \(T\) (see Lemma 3.1 in \cite{GutmanMath}). 
\\

The \emph{Randi\'c (connectivity) index} of a tree \(T\) is defined as \(\sum (deg(i)deg(j))^{-\frac{1}{2}}\), where the summation is over all the edges \(ij\) of \(T\). It was defined by Randi'c \cite{Randic} under the name branching index. Bollob\'as and Erd\H{o}s \cite{BollobasBellaErdos} generalized these indices by defining the \emph{general Randi\'c index} as \(\sum (deg(i)deg(j))^{\alpha}\), where \(\alpha\) is any real number and the summation is again over all the edges \(ij\) of \(T\). Note that for \(\alpha = 1\) this becomes the second Zagreb index, defined above.  
\\

The \emph{Wiener polynomial} (called sometimes Hosoya polynomial or Wiener-Hosoya polynomial) of a tree \(T\) is defined as \(\sum x^{d(i,j)}\), where the sum is taken over all unordered pairs of distinct vertices \((i,j)\) of \(T\)
\cite{Hosoya}, \cite{SaganYehZhang}.

For a tree \(T\), the \emph{partial Wiener polynomial} with respect to a vertex \(w\) of \(T\) is defined as \(\sum x^{d(i,w)}\), where the sum is taken over all vertices \(i\) of \(T\), distinct from \(w\) \cite{DoslicTomsilav}.

\section{Statements and Proofs}

The proofs of our statements are based on examining how a particular statistic S on each of the two trees of Fig. 3 is obtained from the values of S and possibly of some auxiliary statistics on the branches \(\tau(t), \tau(r), \tau(s)\). 
We will use frequently the fact that the degree of the root of \(\tau(n)\) is equal to the number of prime divisors of \(n\), counted with multiplicities.

The sequences obtained by all of the following propositions can be found in OEIS \cite{Encyclopedia} under the indicated number Axxxxxx, where also Maple programs are provided.  

Abusing notation, given a statistic \(S\) on rooted trees, we shall denote \(S(\tau(n))\) by \(S(n)\). In other words, \(S(n)\) is the value of the statistic \(S\) on the rooted tree having Matula number \(n\).   

We would like to point out that some of the statistics we consider have been considered previously in \cite{ElkGutman}, \cite{GutmanMath} and \cite{GutmanIvic}. We include them here in order that the paper be self-contained (some of them play an auxiliary role at other statistics) and because our approach is slightly different. 

The symbols for the various statistics have been selected in the hope that they have at least some limited mnemonic value.  

\begin{prop}
Let \(V\) denote "number of vertices". Then    
\begin{equation}
V(n)=
\begin{cases}
1, &\text{if $n=1$;} \\
1+V(t), &\text{if $n=p_t$;} \\
V(r)+V(s)-1, &\text{if $n=rs$, $r,s \geq 2$.}
\end{cases}
\end{equation}
\end{prop}
\begin{proof}
Follows at once by examining Fig. 3.(A061775) 
\end{proof}

\begin{prop}
Let \(E\) denote "number of edges". Then   
\begin{equation}
E(n)=
\begin{cases}
0, &\text{if $n=1$;} \\
1+E(t), &\text{if $n=p_t$;} \\
E(r)+E(s), &\text{if $n=rs$, $r,s \geq 2$.}
\end{cases}
\end{equation} 
\end{prop}
\begin{proof}
Follows at once by examining Fig. 3. Obviously, \(E(n)=V(n)-1\).  (A196050)
\end{proof}

\begin{prop}
Let \(H\) denote "height". Then   

\begin{equation}
H(n)=
\begin{cases}
0, &\text{if $n=1$;} \\
1+H(t), &\text{if $n=p_t$;} \\
max(H(r),H(s)), &\text{if $n=rs$, $r,s \geq 2$.}
\end{cases}
\end{equation}
\end{prop}
\begin{proof}
Follows at once by examining Fig. 3.   (A109082)
\end{proof}

\begin{prop}
Let \(LLL\) denote "level of lowest leaf". Then   
\begin{equation}
LLL(n)=
\begin{cases}
1+LLL(t), &\text{if $n=p_t$;} \\
min(LLL(r),LLL(s)), &\text{if $n=rs$, $r,s \geq 2$.}
\end{cases}
\end{equation}\end{prop} 

\begin{proof}
Follows at once by examining Fig. 3. The 1-vertex tree, corresponding to \(n=1\), has no leaves.  (A184166)
\end{proof}

\begin{prop}
Let \(LV\) denote "number of leaves". Then   
\begin{equation}
LV(n)=
\begin{cases}
0, &\text{if $n=1$;} \\
1, &\text{if $n=2$;} \\
LV(t), &\text{if $n=p_t$, $t\geq 2$;} \\
LV(r)+LV(s), &\text{if $n=rs$, $r,s \geq 2$.}
\end{cases}
\end{equation}\end{prop} 

\begin{proof}
Follows at once by examining Fig. 3.(A109129)
\end{proof}

\begin{prop}Let \(MD\) denote "maximum vertex degree". Then   

\begin{equation}
MD(n)=
\begin{cases}
0, &\text{if $n=1$;} \\
max(MD(t), 1+\Omega(t)), &\text{if $n=p_t$;} \\
max(MD(r),MD(s), \Omega(r)+\Omega(s)), &\text{if $n=rs$, $r,s \geq 2$.}
\end{cases}
\end{equation}
\end{prop} 

\begin{proof}
Follows at once by examining Fig. 3. (A196046)
\end{proof}

\begin{prop}
Let \(DM\) denote "diameter". Then    
\begin{equation}
DM(n)=
\begin{cases}
0, &\text{if $n=1$;} \\
max(DM(t),1+H(t)), &\text{if $n=p_t$;} \\
max(DM(r),DM(s),H(r)+H(s)), &\text{if $n=rs$, $r,s \geq 2$.}
\end{cases}
\end{equation}\end{prop} 

\begin{proof}
 
(i) If the diameter of \(\tau(t)\) is not the diameter of \(\tau(p_t)\),  then clearly the latter is given by the path going from the root of \(\tau(p_t)\) to a leaf of \(\tau(t)\) of maximum height. 
(ii) Similarly, if neither \(DM(r)\) nor \(DM(s)\) is the diameter of the entire tree \(\tau(rs)\), then clearly this diameter is given by a path passing through the root.  (A196058)
\end{proof}

\begin{prop} Let \(PL\) denote "path length". Then   
\begin{equation}
PL(n)=
\begin{cases}
0, &\text{if $n=1$;} \\
PL(t)+V(t), &\text{if $n=p_t$;} \\
PL(r)+PL(s), &\text{if $n=rs$, $r,s \geq 2$.}
\end{cases}
\end{equation}\end{prop} 

\begin{proof} Each of the \(V(t)\) paths that contribute to \(PL(p_t)\) is 1 unit longer than those that make up \(PL(t)\). The case \(n=rs\) follows at once from Fig. 3b.   (A196047)
\end{proof}

\begin{prop} Let \(EPL\) denote "external path length". Then   
\begin{equation}
EPL(n)=
\begin{cases}
0, &\text{if $n=1$;} \\
1, &\text{if $n=2$;} \\
EPL(t)+LV(t), &\text{if $n=p_t$, $t\geq 2$;} \\
EPL(r)+EPL(s), &\text{if $n=rs$, $r,s \geq 2$.}
\end{cases}
\end{equation}\end{prop} 

\begin{proof} Each of the \(LV(t)\) paths that contribute to \(EPL(p_t)\) is 1 unit longer than those that make up \(EPL(t)\). The case \(n=rs\) follows at once from Fig. 3b.   (A196048)
\end{proof}

\begin{prop} Let \(BV\) denote "number of branching vertices". Then, assuming, without loss of generality, that \(r\) is prime, we have  
\begin{equation}
BV(n)=
\begin{cases}
0, &\text{if $n=1$;} \\
BV(t), &\text{if $n=p_t$, $\Omega(t) \ne 2$;} \\
1+BV(t), &\text{if $n=p_t$, $\Omega(t)=2$;} \\
BV(r)+BV(s), &\text{if $n=rs$, $r,s \geq 2$, $\Omega(s)\ne 2$;} \\
BV(r)+BV(s)+1, &\text{if $n=rs$, $r,s \geq 2$, $\Omega(s)=2$.}  
\end{cases}
\end{equation}\end{prop} 

\begin{proof} (i) In the case \(n=p_t\) we gain a branching vertex only if the root of \(\tau(t)\) has degree 2 in \(\tau(t)\). (ii) In the case \(n=rs\), by joining the trees \(\tau(r)\) and \(\tau(s)\) at their roots, we create a new branching vertex only if the degree of the root of \(\tau(s)\) in the tree \(\tau(s)\) is 2, i.e. \(\Omega(s)=2\).   (A196049)
\end{proof}

\begin{prop} Let \(PV\) denote "number of pendant vertices". Then 
\begin{equation}
PV(n)=
\begin{cases}
0, &\text{if $n=1$;} \\
2, &\text{if $n=2$;} \\
1+LV(t), &\text{if $n=p_t$, $n > 2$;} \\
LV(r)+ LV(s), &\text{if $n=rs$, $r,s \geq 2$. }
\end{cases}
\end{equation}\end{prop} 
\begin{proof} (i) If \(n=p_t\), \(n>2\), then the pendant vertices of \(\tau(n)\) consist of the root and the leaves of \(\tau(t)\). (ii) If \(n\) is composite, \(n=rs\), then the pendant vertices of \(\tau(n)\) consist of the leaves of \(\tau(r)\) and those of \(\tau(s)\).   (A196067)
\end{proof}

\begin{prop} Let \(SP\) denote "number of sibling pairs". Then   
\begin{equation}
SP(n)=
\begin{cases}
0, &\text{if $n=1$;} \\
SP(t), &\text{if $n=p_t$;} \\
SP(r)+SP(s)+\Omega(r)\Omega(s), &\text{if $n=rs$, $r,s \geq 2$.}
\end{cases}
\end{equation}\end{prop} 
\begin{proof} (i) The trees \(\tau(p_t)\) and \(\tau(t)\) have the same number of sibling pairs. (ii) In addition to the sibling pairs of \(\tau(r)\) and those of \(\tau(s)\) we have to count the pairs \((u,v)\), where \(u\) (\(v\)) is a child of the root of \(\tau(rs)\) in \(\tau(r)\) (\(\tau(s)\)).   (A196057)
\end{proof}

\begin{prop} Let \(VL\) denote "visitation length". Then   
\begin{equation}
VL(n)=
\begin{cases}
1, &\text{if $n=1$;} \\
VL(t)+V(t)+1, &\text{if $n=p_t$;} \\
VL(r)+VL(s)-1, &\text{if $n=rs$, $r,s \geq 2$.}
\end{cases}
\end{equation}\end{prop} 

\begin{proof} (i) When going from \(\tau(t)\) to \(\tau(p_t)\) the path length increases by \(V(t)\) and the number of nodes increases by 1.
(ii) The term \(-1\) is due to the fact that in \(VL(r)+VL(s)\) the root is counted twice.   (A196068)
\end{proof}

\begin{prop} Let \(RST\) denote "number of root subtrees". Then   
\begin{equation}
RST(n)=
\begin{cases}
1, &\text{if $n=1$;} \\
1+RST(t), &\text{if $n=p_t$;} \\
RST(r)RST(s), &\text{if $n=rs$, $r,s \geq 2$.}
\end{cases}
\end{equation}\end{prop} 

\begin{proof} (i) The root subtrees of \(\tau(p_t)\) are the "extended" root subtrees of \(\tau(t)\) and the root of \(\tau(p_t)\) (a 1-vertex tree).
(ii) Each root subtree of \(\tau(rs)\) is obtained by joining a root subtree of \(\tau(r)\) and a root subtree of \(\tau(s)\) at their roots. (A184160)
\end{proof}

\begin{prop} Let \(ST\) denote "number of subtrees". Then   
\begin{equation}
ST(n)=
\begin{cases}
1, &\text{if $n=1$;} \\
1+ST(t)+RST(t), &\text{if $n=p_t$;} \\
ST(r)+ST(s)+(RST(r)-1)(RST(s)-1)-1, &\text{if $n=rs$, $r,s \geq 2$.}
\end{cases}
\end{equation}\end{prop} 

\begin{proof} (i) As shown in the proof of the previous proposition, \(1+RST(t)\) is the number of  subtrees of \(\tau(p_t)\) that do contain the root; those that do not contain the root are clearly the subtrees of \(\tau(t)\), counted by \(ST(t)\). (ii) The subtrees that belong entirely to \(\tau(r)\) or entirely to \(\tau(s)\) are counted by \(ST(r)+ST(s)-1\); the remaining ones are obtained by joining at their roots a non-1-vertex root subtree of \(\tau(r)\) with a non-1-vertex root subtree of \(\tau(s)\). (A84161)
\end{proof}

The statistic, depending on a real parameter \(\alpha\), considered in the following proposition is needed at the statistics "2nd Zagreb index" and "general Randi\'c index".

\begin{prop} Let \(A_{\alpha} = \sum deg(i)^{\alpha}\), where the summation extends over all vertices at level 1. Then  
\begin{equation}
A_{\alpha}(n)=
\begin{cases}
0, &\text{if $n=1$;} \\
(1+\Omega(t))^{\alpha}, &\text{if $n=p_t$;} \\
A_{\alpha}(r)+A_{\alpha}(s), &\text{if $n=rs$, $r,s \geq 2$.}
\end{cases}
\end{equation}\end{prop} 

\begin{proof} In \(\tau(p_t)\), there is only one vertex at level 1; its degree is \(1+\Omega(t)\). 
The other cases are obviously true.   (A196052)
\end{proof}

\begin{prop} Let \(W\) denote "Wiener index". Then   
\begin{equation}
W(n)=
\begin{cases}
0, &\text{if $n=1$;} \\
W(t)+PL(t)+E(t)+1, &\text{if $n=p_t$;} \\
W(r)+W(s)+PL(r)E(s)+PL(s)E(r), &\text{if $n=rs$, $r,s \geq 2$.}
\end{cases}
\end{equation}\end{prop} 
\begin{proof}
(i) In the case \(n=p_t\), the first term takes care of the distances within \(\tau(t)\) while the last 3 terms give the sum of the distances between the vertices of \(\tau(t)\) and the root. (ii) In the case \(n=rs\) the first two terms take care of the distances within \(\tau(r)\) and those within \(\tau(s)\). In the sum of the distances between the vertices of \(\tau(r)\) and those of \(\tau(s)\), the sum \(PL(r)\) (\(PL(s)\)) occurs for each nonroot vertex of \(\tau(s)\) (\(\tau(r)\)) i.e. \(E(s)\) (\(E(r)\)) times.   (A196051) 

\end{proof}

\begin{prop} Let \(TW\) denote "terminal Wiener index". Then, assuming, without loss of generality, that \(r\) is prime, we have  
\begin{equation}
TW(n)=
\begin{cases}
0, &\text{if $n=1$;} \\
1, &\text{if $n=2$;} \\
TW(t)+LV(t), &\text{if $n=p_t$ and } \\
&\text{$\Omega(t)=1$ i.e. $t$ is prime}\\
TW(t)+EPL(t)+LV(t),  &\text{if $n=p_t$} \\
&\text{and $\Omega(t)>1$;}\\
TW(r)-EPL(r)+TW(s)-& \\
EPL(s)+EPL(r)LV(s)+EPL(s)LV(r), &\text{if $n=rs$, $r,s \geq 2$, $s$ is prime} \\
TW(r)-EPL(r)+TW(s)+& \\
EPL(r)LV(s)+EPL(s)LV(r), &\text{if $n=rs$, $r,s \geq 2$, $s$ is not prime.}
\end{cases}
\end{equation}\end{prop} 

\begin{proof} 
(i) If both \(n\) and \(t\) are prime, then the distance between each leaf of \(\tau(t)\) and the root of \(\tau(t)\) has to be increased by 1 unit, explaining the term \(LV(t)\). (ii) If \(n=p_t\) and \(\Omega(t)>1\), then the pendant vertices of \(\tau(p_t)\) consist of  those of \(\tau(t)\) and the root of \(\tau(p_t)\). Consequently, to \(TW(t)\) we have to add the sum of the additional distances, namely \(EPL(t)+LV(t)\). (iii) In the case \(n=rs\), \(\Omega(s)=1\),  the root of \(\tau(rs)\) is definitely not a pendant vertex. Consequently, the contribution to \(TW(rs)\) of the distances between the pendant vertices within \(\tau(r)\) is \(TW(r)-EPL(r)\); same explanation for \(TW(s)-EPL(s)\). In the sum of the distances between the leaves of \(\tau(r)\) and those of \(\tau(s)\), the sum \(EPL(r)\) occurs \(LV(s)\) times and the sum \(EPL(s)\) occurs \(LV(r)\) times. (iv) The case \(n=rs\), \(\Omega(s)>1\) is similar to the previous one, except that now the root of \(\tau(s)\) is not a pendant vertex in \(\tau(s)\) and, therefore, there is no need to correct \(TW(s)\) by a subtraction.   (A196055) 
 \end{proof}

\begin{prop} Let \(Z1\) denote "first Zagreb index". Then   
\begin{equation}
Z1(n)=
\begin{cases}
0, &\text{if $n=1$;} \\
Z1(t)+2+2\Omega(t), &\text{if $n=p_t$;} \\
Z1(r)+Z1(s)-\Omega(r)^2-\Omega(s)^2+\Omega(n)^2, &\text{if $n=rs$, $r,s \geq 2$.}
\end{cases}
\end{equation}\end{prop} 

\begin{proof}
(i) In the case \(n=p_t\), to \(Z1(t)\) we add 1 (the squared degree of the root of \(\tau(p_t)\) and we take into account that the degree of the root of \(\tau(t)\) has increased by 1 unit (i.e. we add \((1+\Omega(t))^2-\Omega(t)^2\)). (ii) In the case \(n=rs\) we take into account that the degree of the root is \(\Omega(n)\).   (A196053) 
\end{proof}

\begin{prop} Let \(Z2\) denote "second Zagreb index". Then   
\begin{equation}
Z2(n)=
\begin{cases}
0, &\text{if $n=1$;} \\
Z2(t)+A(t)+\Omega(t)+1, &\text{if $n=p_t$;} \\
Z2(r)+Z2(s)+A(r)\Omega(s)+A(s)\Omega(r), &\text{if $n=rs$, $r,s \geq 2$.}
\end{cases}
\end{equation}\end{prop}

\begin{proof} Here we need the auxiliary statistic "sum of degrees of vertices at level 1", denoted by \(A_1\) and considered in proposition 3.16.  
(i) Assume that \(n=p_t\). Let \(u\) be the root of \(\tau(t)\) and let \(ux\) be an edge in 
\(\tau(t)\). This edge contributes \(deg(u)deg(x)\) to \(Z2(t)\) but to \(Z2(p_t)\) it has to 
contribute \(1+deg(u))deg(x)\). The difference is \(deg(x)\), which, summed over all the edges \(ux\) in \(\tau(t)\), yields \(A_1(t)\). The sum \(1+\Omega(t)\) is the contribution of the edge emanating from the root of \(\tau(p_t)\). 
(ii) Assume that \(n=rs\). Let \(v\) be the root of \(\tau(rs)\) and let \(vy\) be an edge in \(\tau(r)\). This edge contributes \(deg_{\tau(r)}(y)deg_{\tau(r)}(v)\) to \(Z2(r)\) but to \(Z2(n)\) it has to contribute \(deg_{\tau(r)}(y)deg_{\tau(n)}(v)\). The difference is \(deg_{\tau(r)}(y)(deg_{\tau(n)}(v)-deg_{\tau(r)}(v))=deg_{\tau(r)}(y)deg_{\tau(s)}(v)=deg_{\tau(r)}(y)\Omega(s)\). Summing over all the edges in \(\tau(r)\) that emanate from \(v\), we obtain \(A_1(r)\Omega(s)\). Same explanation for the edges within \(\tau(s)\).   (A196054) 
\end{proof}

\begin{prop} Let \(NK\) denote "Narumi-Katayama index". Then   
\begin{equation}
NK(n)=
\begin{cases}
0, &\text{if $n=1$;} \\
1, &\text{if $n=2$;} \\
NK(t)(1+\frac{1}{\Omega(t)}), &\text{if $n=p_t$, $t\geq 2$;} \\
NK(r)NK(s)(\frac{1}{\Omega(r)}+\frac{1}{\Omega(s)}), &\text{if $n=rs$, $r,s \geq 2$.}
\end{cases}
\end{equation}\end{prop} 

\begin{proof}
(i) If \(n=p_t\), then in the Narumi-Katayama index \(NK(t)\) of \(\tau(t)\) we replace the degree of the root of \(\tau(t)\) in \(\tau(t)\) by its degree in \(\tau(p_t)\), i.e. \(\Omega(t)\) by \(1+\Omega(t)\). 
(ii) If \(n=rs\), then in the product \(NK(r)NK(s)\) we replace the product \(\Omega(r)\Omega(s)\) of the root degrees in \(\tau(r)\) and \(\tau(s)\) by the root degree \(\Omega(n)=\Omega(r)+\Omega(s)\) in \(\tau(n)\).   (A196063) 
\end{proof}

\begin{prop} Let \(MZ1\) denote "first multiplicative Zagreb index". Then   
\begin{equation}
MZ1(n)=
\begin{cases}
0, &\text{if $n=1$;} \\
1, &\text{if $n=2$;} \\
MZ1(t)(1+\frac{1}{\Omega(t)})^2, &\text{if $n=p_t$, $t\geq 2$;} \\
MZ1(r)MZ1(s)(\frac{1}{\Omega(r)}+\frac{1}{\Omega(s)})^2, &\text{if $n=rs$, $r,s \geq 2$.}
\end{cases}
\end{equation}\end{prop} 

\begin{proof} Since \(MZ1(n)=NK(n)^2\), the proof is entirely similar to that of  proposition 3.21.   (A196065) 
\end{proof}

\begin{prop} Let \(MZ2\) denote "second multiplicative Zagreb index". Then   
\begin{equation}
MZ2(n)=
\begin{cases}
0, &\text{if $n=1$;} \\
1, &\text{if $n=2$;} \\
\frac{MZ2(t)}{\Omega(t)^{\Omega(t)}}(1+\Omega(t))^{1+\Omega(t)}, &\text{if $n=p_t$, $t\geq 2$;} \\
\frac{MZ2(r)MZ2(s)\Omega(n)^{\Omega(n)}}{\Omega(r)^{\Omega(r)}\Omega(s)^{\Omega(s)}}, &\text{if $n=rs$, $r,s \geq 2$.}
\end{cases}
\end{equation}\end{prop} 

\begin{proof}
As mentioned in the introduction, the 2nd multiplicative Zagreb index can be expressed also as 
\(\prod deg(i)^{deg(i)}\), where the product is taken over all vertices \(i\) of \(T\). Making use of this, the proof is similar to the proof of proposition 3.21.   (A196064) 
\end{proof}

\begin{prop} Let \(R_{\alpha}\) denote "general Randi\'c index". Then 
\begin{equation}
R_{\alpha}(n)=
\begin{cases}
0, &\text{if $n=1$;} \\
R_{\alpha}(t)+A_{\alpha}(t)[(1+\Omega(t))^{\alpha}-\Omega(t)^{\alpha}]+ (1+ \Omega(t))^{\alpha}, &\text{if $n=p_t$;} \\
R_{\alpha}(r)+R_{\alpha}(s)+ A_{\alpha}(r)[\Omega(n)^{\alpha}-\Omega(r)^{\alpha}]+ A_{\alpha}(s)[\Omega(n)^{\alpha}-\Omega(s)^{\alpha}] &\text{if $n=rs$, $r,s \geq 2$.}
\end{cases}
\end{equation}\end{prop}

\begin{proof} 
The reasoning is basically the same as that proposition 3.20 which is a special case of this one (\(\alpha = 1\)). (i) In the case \(n=p_t\), the single edge emanating from the root brings in the term \((1+ \Omega(t))^{\alpha}\), while the edges in \(\tau(t)\) have the same contribution to \(R_{\alpha}(n)\) as to \(R_{\alpha}(t)\), except for those emanating from the root of \(\tau(t)\) because for these the degree of the lower endpoints has increased from \(\Omega(t)\) to \(1+\Omega(t)\). (ii) In the case \(n=rs\), the only edges in \(\tau(n)\) which do not have the same contribution to \(R_{\alpha}(n)\) as to  \(R_{\alpha}(r) + R_{\alpha}(s)\) are those emanating from the root. The difference in these contributions is taken care of by the last two terms in the right-hand side.    
\end{proof}

\begin{prop}5. Let \(PWP\) denote "partial Wiener polynomial with respect to the root". Then  
\begin{equation}
PWP(n)=
\begin{cases}
0, &\text{if $n=1$;} \\
x+xPWP(t), &\text{if $n=p_t$;} \\
PWP(r)+PWP(s), &\text{if $n=rs$, $r,s \geq 2$.}
\end{cases}
\end{equation}\end{prop}   

\begin{proof}
If \(n=p_t\), then, in going from \(\tau(t)\) to \(\tau(p_t)\), the distances from the vertices of \(\tau(t)\) to the root of \(\tau(t)\) are increased by 1 unit; this explains the factor \(x\) in \(xPWP(t)\). The additional distance between the root of \(\tau(t)\) and that of \(\tau(p_t)\) explains the term \(x\). The case \(n=rs\) (\(r,s\geq2\)) follows at once from Fig. 3b.   (A196056)    
\end{proof}
\begin{prop} Let \(WP\) denote "Wiener polynomial". Then 
\begin{equation}
WP(n)=
\begin{cases}
0, &\text{if $n=1$;} \\
WP(t)+xPWP(t)+x, &\text{if $n=p_t$;} \\
WP(r)+WP(s)+PWP(r)PWP(s), &\text{if $n=rs$, $r,s \geq 2$.}
\end{cases}
\end{equation}\end{prop} 

\begin{proof}
(i) \(WP(t)\) takes care of the distances within \(\tau(t)\), \(xPWP(t)\) takes care of the distances between the vertices of \(\tau(t)\) and the root of \(\tau(p_t)\), and the term \(x\) is due to the distance between the root of \(\tau(t)\) and the root of \(\tau(p_t)\). 
(ii) In the case \(n=rs\) (\(r,s \geq 2\), the three terms give the contributions of the distances within \(\tau(r)\), within \(\tau(s)\), and between the vertices in \(\tau(r)\) and those in \(\tau(s)\), respectively.   (A196059) 
\end{proof}

\begin{prop} Let \(DSP\) denote "degree sequence polynomial". Then  
\begin{equation}
DSP(n)=
\begin{cases}
1, &\text{if $n=1$;} \\
DSP(t)+x^{\Omega(t)}(x-1)+x, &\text{if $n=p_t$;} \\
DSP(r)+DSP(s)-x^{\Omega(r)}-x^{\Omega(s)}+x^{\Omega(n)}, &\text{if $n=rs$, $r,s \geq 2$.}
\end{cases}
\end{equation}\end{prop} 

\begin{proof}
(i) When we go from \(\tau(t)\) to \(\tau(p_t)\), the degree of the root of \(\tau(t)\) increases by 1 unit, explaining the term \(x^{\Omega(t)}(x-1)=-x^{\Omega(t)}+x^{1+\Omega(t)}\). We also have an additional vertex of degree 1 (the root of \(\tau(p_t)\)), explaining the term \(x\). 
(ii) We are replacing the terms corresponding to the root of \(\tau(r)\) and to the root \(\tau(s)\) by the  term corresponding to the root of \(\tau(n)\).(A182907) 
\end{proof}

\begin{prop} Let \(EDP\) denote "exit-distance polynomial". Then  
\begin{equation}
EDP(n)=
\begin{cases}
1, &\text{if $n=1$;} \\
EDP(t)+x^{1+LLL(t)}, &\text{if $n=p_t$;} \\
EDP(r)+EDP(s)-x^{max(LLL(r),LLL(s))} &\text{if $n=rs$, $r,s \geq 2$.}
\end{cases}
\end{equation}\end{prop} 
\begin{proof}
(i) The vertices in \(\tau(t)\) have the same exit distances in \(\tau(p_t)\) as in \(\tau(t)\) while the root of \(\tau(p_t)\) has exit distance \(1+LLL(t)\); (ii) from the two terms corresponding to the roots of \(\tau(r)\) and \(\tau(s)\) we remove the term corresponding to the larger exit distance. (A184167) 
\end{proof}

{\bf Remark}
From the propositions that consider polynomial-valued statistics we can re-obtain the results regarding several of the above considered statistics. Without going into details,  

(i) denoting \(f(x)=PWP(n)\), we have \(H(n)=\) degree of \(f(x)\), \(E(n)=f(1)\), and \(PL(n)=\frac{df}{dx}|_{x=1}\); also, the number of vertices at level \(k\) \((k\geq 1)\) in the tree \(\tau(n)\) is equal to \([x^k]f(x)\);

(ii) denoting \(g(x)=WP(n)\), we have \(DM(n)=\) degree of \(g(x)\) and \(W(n)= \frac{dg}{dx}|_{x=1}\).

(iii) denoting \(h(x)=DSP(n)\), we have \(V(n)=h(1)\), \(MD(n)=\) degree of \(h(x)\), \(PV(n)=[x]h(x)\), and \(BV(n)=h(1)-[x]h(x)-[x^2]h(x)\).

Moreover, making use of some of these polynomial-valued statistics, we can consider also some statistics that are not included in the propositions given above. 

(a) \emph{The hyper-Wiener index}. It has been proved in \cite{CashGG} that the hyper-Wiener index of a connected graph is equal to \(g'(1)+\frac{1}{2}g"(1)\), where \(g(x)\) is the Wiener polynomial of \(G\) (A196060).

(b) \emph{The multiplicative Wiener index}. From the definition there follows easily that the multiplicative Wiener index of a graph \(G\) is equal to \(\prod_{k>1} k^{\delta(G,k)}\), where 
\(\delta(G,k)\) is the number of vertex pairs in \(G\) that are at distance \(k\) (see Eq. (5a) in \cite{GutmanLinertLukovitsTomovic} (A196061).   

(c) \emph{The Wiener polarity index}. This is the coefficient of \(x^3\) in the Wiener polynomial of the graph (A184156). Obviously, taking other coefficients in the Wiener polynomial, one obtains generalized Wiener polarity indices (see \cite{IlicIlic}).

(d) In \cite{IvanciucIvanciucKleinSeitzBalaban} one introduces a family of statistics on a connected graph \(G\), based on the graph distances. In particular, Sum\(E(1,G)\) (Sum\(O(1,G)\)) denotes the sum of the even (odd) distances between unordered pairs of vertices of \(G\). Clearly, it can be obtained by evaluating at \(x=1\) the derivative of the even (odd) part of the Wiener polynomial of \(G\) (A184157, A184158).

(e) Denoting \(m(x)=EDP(n)\), it is easy to see that (i) the sum of the exit distances of all vertices of \(\tau(n)\) is equal to \(\frac{dm}{dx}|_{x=1}\) (A184168) ,(ii) the maximum exit distance over the vertices of \(\tau(n)\) is the degree of the polynomial \(m(x)\) (A184169), and (iii) the number of vertices of \(\tau(n)\) having maximum exit distance is the coefficient of the highest power of \(x\) in \(m(x)\) (A184170). 
\newpage

\bibliography{matula.bib}
\bibliographystyle{plain}
\end{document}

%% file: fig1.tex
%\documentclass{article}
%\usepackage{pgf}
%\usepackage{xcolor}

%\begin{document}
\begin{pgfpicture}{0cm}{0cm}{10cm}{4cm}
\pgfsetstartarrow{\pgfarrowcircle{3pt}}
\pgfsetendarrow{\pgfarrowcircle{2pt}}
\pgfline{\pgfxy(5,0)}{\pgfxy(2,2)}
\pgfputat
{\pgflabel{.9}{\pgfxy(5,0)}{\pgfxy(2,2)}{3pt}}
{\pgfbox[right,down]{$v_1$}}
\pgfclearstartarrow
\pgfline{\pgfxy(5,0)}{\pgfxy(4.5,2)}
\pgfputat
{\pgflabel{.9}{\pgfxy(5,0)}{\pgfxy(4.5,2)}{3pt}}
{\pgfbox[right,down]{$v_2$}}
\pgfline{\pgfxy(5,0)}{\pgfxy(8,2)}
\pgfputat
{\pgflabel{.9}{\pgfxy(5,0)}{\pgfxy(8,2)}{3pt}}
{\pgfbox[right,down]{$v_d$}}

\pgfputat{\pgfxy(6.3,2)}{\pgfbox[center,center]{\Huge{. . .}}}
\pgfcircle[stroke]{\pgfxy(2,3)}{1cm}
\pgfcircle[stroke]{\pgfxy(4.5,3)}{1cm}
\pgfcircle[stroke]{\pgfxy(8,3)}{1cm}
\pgfputat{\pgfxy(2,3)}{\pgfbox[center,center]{$T_1$}}
\pgfputat{\pgfxy(4.5,3)}{\pgfbox[center,center]{$T_2$}}
\pgfputat{\pgfxy(8,3)}{\pgfbox[center,center]{$T_d$}}
\end{pgfpicture}
%\end{document}

%% file: fig2.tex
%\documentclass{article}
%\usepackage{pgf}
%\usepackage{xcolor}
%\usepackage[margin=1cm]{geometry}
%\begin{document}
\begin{pgfpicture}{0cm}{-1cm}{24cm}{5cm}
\pgfsetstartarrow{\pgfarrowcircle{2pt}}
\pgfsetendarrow{\pgfarrowcircle{2pt}}
\begin{pgfmagnify}{.75}{.75}
\pgfputat{\pgfxy(7,-.1)}{\pgfbox[center,top]{$987654321=3\cdot3\cdot17\cdot17\cdot379721$}}
\pgfline{\pgfxy(7,0)}{\pgfxy(0,1)}
\pgfputat
{\pgflabel{1}{\pgfxy(7,0)}{\pgfxy(0,1)}{3pt}}
{\pgfbox[right,top]{2}}
\pgfclearstartarrow
\pgfline{\pgfxy(0,1)}{\pgfxy(0,2)}
\pgfputat
{\pgflabel{1}{\pgfxy(0,1)}{\pgfxy(0,2)}{3pt}}
{\pgfbox[right,top]{1}}
\pgfline{\pgfxy(7,0)}{\pgfxy(3,1)}
\pgfputat
{\pgflabel{1}{\pgfxy(7,0)}{\pgfxy(3,1)}{3pt}}
{\pgfbox[right,top]{2}}
\pgfline{\pgfxy(7,0)}{\pgfxy(7,1)}
\pgfputat
{\pgflabel{1}{\pgfxy(7,0)}{\pgfxy(7,1)}{3pt}}
{\pgfbox[right,top]{7}}
\pgfline{\pgfxy(7,0)}{\pgfxy(10,1)}
\pgfputat
{\pgflabel{1}{\pgfxy(7,0)}{\pgfxy(10,1)}{-3pt}}
{\pgfbox[left,top]{7}}
\pgfline{\pgfxy(7,0)}{\pgfxy(17,1)}
\pgfputat
{\pgflabel{.95}{\pgfxy(7,0)}{\pgfxy(17,1)}{-6pt}}
{\pgfbox[left,top]{$32277=3\cdot7\cdot29\cdot53$}}
\pgfline{\pgfxy(3,1)}{\pgfxy(3,2)}
\pgfputat
{\pgflabel{1}{\pgfxy(3,1)}{\pgfxy(3,2)}{3pt}}
{\pgfbox[right,top]{1}}
\pgfline{\pgfxy(7,1)}{\pgfxy(7,2)}
\pgfputat
{\pgflabel{1}{\pgfxy(7,1)}{\pgfxy(7,2)}{3pt}}
{\pgfbox[right,top]{$4=2\cdot2$}}
\pgfline{\pgfxy(10,1)}{\pgfxy(10,2)}
\pgfputat
{\pgflabel{1}{\pgfxy(10,1)}{\pgfxy(10,2)}{3pt}}
{\pgfbox[right,top]{$4=2\cdot2$}}
\pgfline{\pgfxy(17,1)}{\pgfxy(12,2)}
\pgfputat
{\pgflabel{1}{\pgfxy(17,1)}{\pgfxy(12,2)}{3pt}}
{\pgfbox[right,top]{2}}
\pgfline{\pgfxy(17,1)}{\pgfxy(15,2)}
\pgfputat
{\pgflabel{.95}{\pgfxy(17,1)}{\pgfxy(15,2)}{-7pt}}
{\pgfbox[left,top]{$4=2\cdot2$}}
\pgfline{\pgfxy(17,1)}{\pgfxy(18,2)}
\pgfputat
{\pgflabel{1}{\pgfxy(17,1)}{\pgfxy(18,2)}{-3pt}}
{\pgfbox[left,center]{$10=2\cdot 5$}}
\pgfline{\pgfxy(17,1)}{\pgfxy(22,2)}
\pgfputat
{\pgflabel{1}{\pgfxy(17,1)}{\pgfxy(22,2)}{-4pt}}
{\pgfbox[left,top]{$16=2\cdot2\cdot2\cdot2$}}
\pgfline{\pgfxy(7,2)}{\pgfxy(6,3)}
\pgfputat
{\pgflabel{1}{\pgfxy(7,2)}{\pgfxy(6,3)}{3pt}}
{\pgfbox[right,top]{1}}
\pgfline{\pgfxy(7,2)}{\pgfxy(8,3)}
\pgfputat
{\pgflabel{1}{\pgfxy(7,2)}{\pgfxy(8,3)}{-3pt}}
{\pgfbox[left,top]{1}}
\pgfline{\pgfxy(10,2)}{\pgfxy(9,3)}
\pgfputat
{\pgflabel{1}{\pgfxy(10,2)}{\pgfxy(9,3)}{3pt}}
{\pgfbox[right,top]{1}}
\pgfline{\pgfxy(10,2)}{\pgfxy(11,3)}
\pgfputat
{\pgflabel{1}{\pgfxy(10,2)}{\pgfxy(11,3)}{-3pt}}
{\pgfbox[left,top]{1}}
\pgfline{\pgfxy(12,2)}{\pgfxy(12,3)}
\pgfputat
{\pgflabel{1}{\pgfxy(12,2)}{\pgfxy(12,3)}{-3pt}}
{\pgfbox[left,top]{1}}
\pgfline{\pgfxy(15,2)}{\pgfxy(14,3)}
\pgfputat
{\pgflabel{1}{\pgfxy(15,2)}{\pgfxy(14,3)}{3pt}}
{\pgfbox[right,top]{1}}
\pgfline{\pgfxy(15,2)}{\pgfxy(16,3)}
\pgfputat
{\pgflabel{1}{\pgfxy(15,2)}{\pgfxy(16,3)}{-3pt}}
{\pgfbox[left,top]{1}}
\pgfline{\pgfxy(18,2)}{\pgfxy(17,3)}
\pgfputat
{\pgflabel{1}{\pgfxy(18,2)}{\pgfxy(17,3)}{3pt}}
{\pgfbox[right,top]{1}}
\pgfline{\pgfxy(18,2)}{\pgfxy(19,3)}
\pgfputat
{\pgflabel{1}{\pgfxy(18,2)}{\pgfxy(19,3)}{-3pt}}
{\pgfbox[left,top]{3}}
\pgfline{\pgfxy(19,3)}{\pgfxy(19,4)}
\pgfputat
{\pgflabel{1}{\pgfxy(19,3)}{\pgfxy(19,4)}{-3pt}}
{\pgfbox[left,top]{2}}
\pgfline{\pgfxy(19,4)}{\pgfxy(19,5)}
\pgfputat
{\pgflabel{1}{\pgfxy(19,4)}{\pgfxy(19,5)}{-3pt}}
{\pgfbox[left,top]{1}}
\pgfline{\pgfxy(22,2)}{\pgfxy(20.5,3)}
\pgfputat
{\pgflabel{1}{\pgfxy(22,2)}{\pgfxy(20.5,3)}{3pt}}
{\pgfbox[right,top]{1}}
\pgfline{\pgfxy(22,2)}{\pgfxy(21.5,3)}
\pgfputat
{\pgflabel{1}{\pgfxy(22,2)}{\pgfxy(21.5,3)}{3pt}}
{\pgfbox[right,top]{1}}
\pgfline{\pgfxy(22,2)}{\pgfxy(22.5,3)}
\pgfputat
{\pgflabel{1}{\pgfxy(22,2)}{\pgfxy(22.5,3)}{-3pt}}
{\pgfbox[left,top]{1}}
\pgfline{\pgfxy(22,2)}{\pgfxy(23.5,3)}
\pgfputat
{\pgflabel{1}{\pgfxy(22,2)}{\pgfxy(23.5,3)}{-3pt}}
{\pgfbox[left,top]{1}}

\end{pgfmagnify}

\end{pgfpicture}
%\end{document}

%% file: fig3.tex
%\documentclass{article}
%\usepackage{pgf}
%\usepackage{xcolor}

%\begin{document}
\begin{pgfpicture}{0cm}{0cm}{10cm}{4cm}
\pgfputat{\pgfxy(0,0)}{\pgfbox[center,center]{$a)$}}
\pgfputat{\pgfxy(2,1)}{\pgfbox[center,top]{$\tau(p_t)$}}
\pgfputat{\pgfxy(7.5,0)}{\pgfbox[center,center]{$b)$}}
\pgfputat{\pgfxy(9,1)}{\pgfbox[center,top]{$\tau(rs)$}}
\pgfline{\pgfxy(2,1)}{\pgfxy(2,2)}
\pgfcircle[stroke]{\pgfxy(2,3)}{1cm}
\pgfputat{\pgfxy(2,3)}{\pgfbox[center,center]{$\tau(t)$}}

\begin{pgftranslate}{\pgfxy(1.8,-4.1)}
\begin{pgfrotateby}{\pgfdegree{30}}
\pgfmoveto{\pgfxy(9,1)}
\pgfcurveto{\pgfxy(7,2.5)}{\pgfxy(9,4)}{\pgfxy(9,3)}
\pgfclosepath
\pgfstroke
\end{pgfrotateby}
\end{pgftranslate}

\begin{pgftranslate}{\pgfxy(.75,4.9)}
\begin{pgfrotateby}{\pgfdegree{-30}}
\pgfmoveto{\pgfxy(9,1)}
\pgfcurveto{\pgfxy(11,2.5)}{\pgfxy(9,4)}{\pgfxy(9,3)}
\pgfclosepath
\pgfstroke
\end{pgfrotateby}
\end{pgftranslate}
\pgfputat{\pgfxy(8,2.2)}{\pgfbox[center,center]{$\tau(r)$}}
\pgfputat{\pgfxy(10,2.2)}{\pgfbox[center,center]{$\tau(s)$}}
\end{pgfpicture}

%\end{document}